\documentclass[11pt,fleqn]{article}
\usepackage{graphicx}
\usepackage{fullpage}
\usepackage{epsfig}
\usepackage{color}
\usepackage{latexsym}
\usepackage{amssymb}
\usepackage{amsmath}
\usepackage[new]{old-arrows}
\usepackage{xspace}
\usepackage[mathscr]{euscript}
\usepackage{cite}
\usepackage{url}
\usepackage{latexsym,epsfig,subfigure}
\usepackage{color}
\usepackage{changepage}
\usepackage{tabto}
\usepackage{bm}
\usepackage[shortlabels]{enumitem}
\DeclareSymbolFont{msbm}{U}{msb}{m}{n}
\DeclareMathSymbol{\C}{\mathalpha}{msbm}{'103}
\DeclareMathSymbol{\R}{\mathalpha}{msbm}{'122}
\DeclareMathSymbol{\Q}{\mathalpha}{msbm}{'121}
\DeclareMathSymbol{\Z}{\mathalpha}{msbm}{'132}
\DeclareMathSymbol{\N}{\mathalpha}{msbm}{'116}
\DeclareMathSymbol{\K}{\mathalpha}{msbm}{'113}

\newcommand{\old}[1]{{}}

\newcommand{\bi}{\begin{itemize}}
\newcommand{\ei}{\end  {itemize}}
\newcommand{\bt}{\begin{tabbing}}
\newcommand{\et}{\end  {tabbing}}
\newcommand{\be}{\begin{enumerate}}
\newcommand{\ee}{\end  {enumerate}}

\makeatletter
\def\begin@lgo{\begin{minipage}{1in}\begin{tabbing}
        \quad\=\qquad\=\qquad\=\qquad\=\qquad\=\qquad\=\qquad\=\kill}
\def\end@lgo{\end{tabbing}\end{minipage}}

\makeatother

\makeatletter
\@ifundefined{abovecaptionskip}{\newlength\abovecaptionskip}
\long\def\@makecaption#1#2{
   \vskip \abovecaptionskip
   \setbox\@tempboxa\hbox{{\sf\footnotesize \textbf{#1.} #2}}
   \ifdim \wd\@tempboxa >\hsize         
       {\sf\footnotesize \textbf{#1.} #2\par}
     \else                              
       \hbox to\hsize{\hfil\box\@tempboxa\hfil}
   \fi}
\makeatother

\title{A Finitist's Manifesto: Do we need to Reformulate \\the Foundations of Mathematics?}
\author{Jonathan Lenchner}

\date{}

\begin{document}

\maketitle

\section{The Problem}
There is a problem with the foundations of classical mathematics, and potentially even with the foundations of computer science, that mathematicians have by-and-large ignored. This essay is a call for practicing mathematicians who have been sleep-walking in their infinitary mathematical paradise to take heed. 
Much of mathematics relies upon either (i) the ``existence'' of objects that contain an infinite number of elements, (ii) our ability, ``in theory'', to compute with an arbitrary level of precision, or (iii) our ability, ``in theory'', to compute for an arbitrarily large number of time steps.  All of calculus relies on the notion of a limit. The monumental results of real and complex analysis rely on a seamless notion of the ``continuum'' of real numbers, which extends in the plane to the complex numbers and gives us, among other things, ``rigorous'' definitions of continuity, the derivative, various different integrals, as well as the fundamental theorems of calculus and of algebra -- the former of which says that the derivative and integral can be viewed as inverse operations, and the latter of which says that every polynomial over $\mathbb{C}$ has a complex root. This essay is an inquiry into whether there is any way to assign meaning to the notions of ``existence'' and ``in theory'' in (i) to (iii) above. 

On the one hand, we know from quantum mechanics that making arbitrarily precise measurements of objects is impossible. By the Heisenberg Uncertainty Principle the moment we pin down an object, typically an elementary particle, in space, thereby bringing its velocity, and hence momentum, down to near $0$, there is a limit to how precisely we can measure its spatial coordinates. In symbols:
\begin{equation*}
	\sigma_x \sigma_p \geq \frac{h}{4\pi},
\end{equation*}
where $\sigma_x$ is the standard deviation of position, $\sigma_p$ is the standard deviation of momentum, and $h$ is Planck's constant. Equally troubling, though not as widely noted, is that we believe the universe is finite, so if we could turn the universe into a vast computing device, or even under the assumption that it \textit{is} such a device, there would only be a finite amount of space to write out the decimal places of any measurement. 

I claim that physics presents us with an epistemological problem\footnote{Epistemological problems in philosophy ask about the nature of our ``knowing'' and how we know what we know. See, for example, the discussion in \cite{POLLOCK:1968}.} with our conception of the foundations of mathematics. 
Suppose I challenge you and say that despite everything you know, there is actually a largest natural number (in other words, a number of the form $0, 1,...$), though I don't quite have the facilities to show it to you. How do you prove to me that I am wrong? You would like to say, ``take your $N$, form $N+1$, and it will be bigger.'' However, in return I tell you that to deliver $N$ to you, it would take up all the storage in the universe, and, moreover, we would know that I am representing $N$ in the most compact format possible (e.g., perhaps I am using iterated exponentiation, or even that I am using some of my space to define new symbols that are shorthand for some amount of repeated exponentiation and other things, and using all this to define $N$). You have no ability to even talk to me about $N+1$. Does this number ``exist'' or not? You push me further and say that you are talking about the number $N+1$ in a ``conceptual'' way; after all -- every number has a successor -- it is one of the postulates of Peano Arithmetic! But what do you mean by such a conceptual way? For this $N$ there is no mathematical way whatsoever to talk about $N+1$ (or any number larger than $N$ for that matter). So the only way to talk about $N+1$ is informally -- through our informal language. But even allowing informal language, given our finite vocabularies and space limitations there are only finitely many natural numbers we can describe in this way, and again a maximum such number. While we might be able to get away with a phrase of the form ``the largest natural number expressible using formal mathematics plus one'' allowing the phrase ``the largest number we are able to express (informally) plus one'' would be obvious rubbish. Thus, perhaps we must have a gentlemen's agreement never to discuss anything beyond this largest informally expressible number.  Whether such a resolution is satisfactory (especially to gentlewomen) I leave for the reader to decide. 

Let us now, however, circle back to formal mathematics. 
After all, an expression of the form $\displaystyle{\lim_{n \rightarrow \infty}} f(n)$, for some function $f$, is a formal mathematical expression, not something informal! You may say that as we take larger and larger $n$ we are talking about going from $n$ to $n+1$ in a \textit{general} sense, in the sense of $n$ being \textit{any} natural number, not in a \textit{particular} natural number, let alone the case of our special number $N$. For example, you say, I really don't need to talk about $N+1$ in this particular instance since in all arguments about limits I am using some form of implicit or explicit induction and the $n$ in these arguments does not refer to specific values of $n$.  For example, in the explicit induction proving that 
\begin{equation} \label{eqn:gauss}
	\sum_{i=1}^K i = \frac{K(K+1)}{2},
\end{equation}
we establish the result for a \textit{generic} $n$ and then prove that it holds also for the \textit{generic} $n+1$. Thus, first stopping to observe that (\ref{eqn:gauss}) holds for $K = 0$, we then assume it holds for $K = n$, in other words,
\begin{equation*}
	\sum_{i=1}^n i = \frac{n(n+1)}{2}.
\end{equation*}
Then for $K = n+1$ we establish
\begin{eqnarray*}
	\sum_{i=1}^{n+1} i & = & \sum_{i=1}^n i + (n+1) \\
	& = & \frac{n(n+1)}{2} + \frac{2(n+1)}{2} \\
	& = & \frac{(n+2)(n+1)}{2}.
\end{eqnarray*}
And we have just used $n+1$ (and $n+2$) in this ``generic'' sense.  
To this argument I say yes, what you have shown is indeed that the familiar expression (\ref{eqn:gauss}) holds, at least for all \textit{practical} values of $K$, and my intention is not to challenge all of existing finitary mathematics. However, let us now return to the expression $\displaystyle{\lim_{n \rightarrow \infty}} f(n)$ and put $L = \displaystyle{\lim_{n \rightarrow \infty}} f(n)$. To a mathematician, such an expression is equivalent to saying that 
for any arbitrarily small epsilon there is a point beyond which $f(n)$ is within epsilon of the limiting value $L$. In formal mathematics this notion is expressed thusly: $(\forall \epsilon > 0) \exists m (\forall n>m)~|L - f(n)| < \epsilon$. Hence, the mathematician argues, for any given $k$ one may throw away the first $k$ values and the limiting value, $L$, must be precisely the same. But what if we throw away everything up until our problematic $N$? We end up with an induction that cannot get started! It is as if we have (\ref{eqn:gauss}) but cannot even check the case $K=0$.

Lest you think this is a problem just of the arbitrarily large, note too that we have exactly the same problem when talking about the arbitrarily small, for example, when taking any derivative, or trying to evaluate any other expression of the form  
$\displaystyle{\lim_{\epsilon \rightarrow 0} f(\epsilon)}$.
The only difference in the case of the arbitrarily small is that we have to store, or otherwise represent, a potentially unlimited number of digits to the right of the decimal point. There is therefore a smallest definable $\epsilon > 0$ just like there is a largest definable $N$. 

\medskip


Physics therefore imposes a limit to the meaningfulness of mathematical assertions and we need to question the sanctity of the edifice upon which much of classical mathematics is built.  Despite the amazing structures that have been built atop it -- constructs like Riemannian and semi-Riemannian manifolds, the backbone of our understanding of general relativity and gravity, Hilbert Spaces of operators, the framework of our understanding of quantum mechanics, and many others, and all the wonderful theorems and further constructs that have been built atop this firmament -- the foundations are in peril.
Many a mathematician will bristle at this conclusion, since the aforementioned constructs speak to them like a great painting does to an art connoisseur, egging him or her on with the beauty of the subject, as if to say that something so beautiful cannot \textit{possibly} be built on an illusory edifice. 

\medskip

If we leave our culturally acquired framework for understanding infinity and arbitrarily large numbers at the doorstep, can we define these notions from more primitive ones in a way that will satisfy the extreme finitist -- someone who denies the ability to imagine arbitrarily large space or arbitrarily long time periods (and for whom an attempt to establish the former notions by assuming these latter two seems circular)?

\section{A Prior Crisis at the Foundations}

The dilemma we now face is not dissimilar to that mathematicians had to confront at the turn of the 20th century, when Bertrand Russell unveiled his famous paradox that challenged the then prevalent informal use of set theory.  Prior to Russell's discovery, mathematicians used the notion of a set completely informally, and one could form collections of anything that mathematicians could think of and call the result a ``set.'' Russell asked us to consider the set $\mathscr{N}$ of all sets that do not contain themselves (the set of all so-called ``normal'' sets).  He then further asked whether the set $\mathscr{N}$ contained itself or not. If $\mathscr{N} \in \mathscr{N}$ then $\mathscr{N}$ contains a set that contains itself, which is contrary to its definition. On the other hand if  $\mathscr{N} \notin \mathscr{N}$ then $\mathscr{N}$ does not contain all normal sets since $\mathscr{N}$ is one of them. One is thus led to a contradiction either way. There must, therefore, have been something wrong with the definition of $\mathscr{N}$ in the first place -- and hence to our informal notion of what it means to \textit{be} a ``set.''

To get a sense of just how serious and unsettling Russell's paradox was for mathematicians of the day, note that its revelation came in 1901, just a year after Hilbert's famous address at the International Congress of Mathematicians where he unveiled his celebrated 23 problems, and proclaimed:

\smallskip

\begin{adjustwidth}{2.5em}{2.5em}
``Take any definite unsolved problem, such as the question as to the irrationality of the Euler-Mascheroni constant, or the existence of an infinite number of prime numbers of the form $2^n + 1$. However unapproachable these problems may seem to us, and however helpless we stand before them, we have, nevertheless, the firm conviction that their solution must follow by a finite number of purely logical processes. ...This conviction of the solvability of every mathematical problem is a powerful incentive to the worker. We hear within us the perpetual call: There is the problem. Seek its solution. You can find it by pure reason, for in mathematics there is no \textit{ignorabimus}''\footnote{By the word \textit{ignorabimus}, Hilbert was making implicit reference to the Latin phrase \textit{ignoramus et ignorabimus}, meaning ``we do not know and will not know,'' referring to the fact that scientific knowledge is inherently limited -- an idea popularized by 19th century German physiologist Emil du Bois-Reymond.} \cite{FRANZEN:2005}.
\end{adjustwidth}

\medskip

In response to the Russell paradox, mathematicians created various axiomatizations of Set Theory - rules by which one could form new sets out of existing ones, beginning with a first set, the set with no elements, called the empty set.  These axioms did not provide a means for creating a set of the type $\mathscr{N}$ encountered in Russell's paradox, hence skirting that problem. The most popular axiom system for Set Theory is today known as Zermelo-Fraenkel set theory with the Axiom of Choice, a.k.a, ZFC \cite{JECH:1997}. But there is nothing completely ``obvious'' about the ZFC axioms, and they have been hotly debated over time. Questions about the Axiom of Choice are legion, but there is also the Axiom of Infinity that in an indirect fashion asserts the existence of infinite sets, and the so-called Axiom Schema of Replacement -- an infinite collection of axioms saying that given any set $X$, and a function $f$ defined on that set, then the image $f(X)$ is also a set. Interestingly, the Axiom Schema of Replacement was not accepted by ZFC co-creator Ernst Zermelo \cite{wiki:replacement} because it gives rise to the Downward L{\"o}wenheim–Skolem theorem of Model Theory, and the paradoxical result that there are countable (set) models of Set Theory. Of further note, Zermelo's contributions to ZFC were his attempt to reduce the so-called Well Ordering Principle, first conceived of by Cantor as a way to prove that any two sets have comparable cardinalities, to simpler principles -- and, in particular, to the customary statement of the Axiom of Choice\footnote{That given any collection of non-empty sets, one can form a new set by picking exactly one element from each of the sets.} -- an axiom that Zermelo claimed to be obvious \cite[p.~46]{TARSKI:2004}.

With the question of what it might mean to be a \textit{set} in the back of our minds, let us circle back and ask a fundamental question: what is a real number? We have a sense of what a natural number is. 
Say we are considering the number 3. We think of ``3'' as an abstraction for various collections of objects that can be put in one-to-one correspondence with one another and have the given number of distinguishable elements. But a real number -- is it an abstraction for any \textit{thing} or just shorthand for some sort of approximate measurement that is getting more and more precise (or perhaps not) as we go, or an expansion of a well defined constant like $\pi$ that, again, is getting more and more precise as we go? Is there any \textit{reality} in the real number itself or is it just a way of summarizing a measurement or approximation \textit{activity}? Does mathematics in this sense blend irrevocably into physics? 

Mathematicians have almost uniformly dodged these questions, believing in the objective reality of mathematical concepts, like the different types of numbers and even the collection of \textit{all} such numbers. Furthermore, in an effort to reduce everything in mathematics to one common foundation, following a precedent set by John von Neumann \cite{vonNeumann:1923},
the natural numbers and then the rationals and the reals are typically created from sets. The number $0$ is represented by the empty set $\emptyset = \{\}$. Then one proceeds to construct $1 = \{0\} = \{\emptyset\}, 2 = \{0,1\} = \{\emptyset, \{\emptyset\}\}$, and so on, ultimately creating the natural numbers, $\mathbb{N}$, the rationals $\mathbb{Q}$, the reals $\mathbb{R}$, the associated arithmetical operations, and beyond -- all simply by employing the axioms of ZFC. Whether this reduction of mathematics to Set Theory is natural or not is an open philosophical question. We have some intuition that various collections of three objects should be considered bona fide sets -- but the number $3$ itself? Are we not reducing a concept that we have ironclad clarity about (i.e., the number $3$) with something (a set) that we have just a very loose intuitive notion of?  
Ironically, the more intuitive notion most of us have about the number `3' being a conceptual shorthand for the collection of all sets with three elements in them -- a notion that can be made non-circular by depicting a single set with three elements and then declaring `3' to be the equivalence class of all sets equinumerous with the given set -- would lead us into trouble, at least in ZFC. If there were a set of all three element sets we would easily be able to create a Russell-type paradox based on it.

%
Before delving further into this subject, it is worth noting that even the most fundamental building block of our contemporary theoretical model of computation makes use of an infinite structure -- the Turing Machine with its infinite tape. The famous P=NP question is a question about limits and does not make sense without our ability to consider arbitrarily large integers. We cannot distinguish P, the class of problems solvable in polynomial time, from EXP, the class of problems solvable in exponential time, in just a finite amount of space. Similarly for the Halting Problem -- all but the most trivial diagonalization arguments use at least countably infinite sets, as does finite model theory with its use of the Compactness Theorem of first order logic, i.e., to show that there is no first order reachability predicate in graph theory\footnote{A result that can be extended to finite graph theory using Ehrenfeuch-Fra\"{i}ss\'{e} games of arbitrary large sizes \cite{LIBKIN:2004} -- in this case at least not requiring any sort of completed infinite set.}. There is no question that some of these results can be salvaged. For example, we can ask whether it is possible to write a program in $n$ bits of space that can tell whether an arbitrary program, also written using $n$ bits of space, will halt. For $n$ of ``modest'' size the answer is \textit{no} and the argument is roughly analogous to the classical relativized diagonalization argument used to prove the Space Hierarchy Theorem \cite{PAPA:1995}. However, when $n$ gets close to the largest expressible natural number, the argument cannot be carried out so it is unclear whether we should accept the argument for ``generic'' $n$ or not.

\section{Intuitionism -- A Better Foundation?} \label{sec:intuitionism}

Recently, the Swiss quantum physicist Nicolas Gisin, most well known for his experimental validation of Bell's Inequality and quantum teleportation, has expressed concerns, much as I have articulated them, in a series of papers \cite{gisin2018indeterminism, GISIN:2019, GISIN:2020} connecting problems with the foundations of mathematics with our conception of time and a constructive mathematical philosophy known as Intuitionism \cite{sep-intuitionism}. A popular account of Gisin's work in this area can be found in \cite{DOES_TIME_REALLY_FLOW:2020}. Intuitionism is a philosophy of mathematics that is in striking contrast to the Platonic Realism view adopted by most mathematicians \cite{SNAPPER:1979a, SNAPPER:1979b}.  
In the classical mathematical view, also known as mathematical formalism, mathematical facts, or theorems, exist and just await our discovery. In the Intuitionistic view, on the other hand, mathematics is entirely a construct of the human mind. The Dutch mathematician  L.E.J. 
Brouwer, the founder of intuitionism and also a great classical mathematician,
initially developed the subject from 1905-1910, thus predating all of contemporary computer science. 
Surely, however, there was nothing special in Brouwer's thinking about the thought process of homo sapiens, so the thinking of other beings, or the in-memory constructs of computational machines, 
would also have been allowed had he done his work in the modern era.  However, the distinction between mathematics as facts to be discovered, and facts that can't possibly exist until they are discovered by a mind, human or otherwise, is a very important one. On the one hand, this notion puts great emphasis on time.  At the turn of the 20th century there was a great debate between Hilbert and Brouwer about the nature of the real numbers \cite{GISIN:2020}. Hilbert promoted the notion that every real number, with its infinite sequence of digits, was a completed object, while Brouwer argued that real numbers just represented a never-ending process that develops over \textit{time} in what he referred to as a ``choice sequence.'' Although Brouwer was supported in his view by the famous mathematician, physicist and philosopher Hermann Weyl \cite{WEYL:1994}, and the great logician Kurt G\"{o}del \cite{GODEL:2003}, Hilbert clearly won this debate and trained mathematicians of today are handed down a Platonic Realist view of real numbers.

Gisin argues, as I have, that because a finite volume of space can only contain a finite amount of information, a real number must also only contain a finite amount of information. He suggests that the notion of a real number be replaced by what he calls a Finite Information Quantity \cite{GISIN:2019}.  Suppose we are expanding a number $x$, between $0$ and $1$, and that we would ordinarily have a base $2$ expansion of the form 
\[
	x = b_1 b_2  ...,
\]
where each $b_i$ is a binary digit -- a $0$ or a $1$. Like a Brouwer choice sequence, these values represent a process that develops over time. However, the numbers $b_i$, for Gisin, are not bits in the conventional sense, but each is a rational number in the closed interval $[0,1]$ giving the probability that the underlying value is a $1$. An underlying probability of $0$ coincides with $b_i$ being definitively $0$. When the values are completely random, so that $b_i = 0.5$, they contain no information. Once these pseudo-bits are pinned down to a $0$ or a $1$ they contain the most information.  The information contained in an entire Finite Information Quantity of this sort is given by\footnote{Where $h(0) = h(1) = 0$ by taking limits, in other words $\lim_{x \rightarrow 0}x \lg{x} = 0$.}
\begin{eqnarray} \label{eqn_info}
	I(x) & = & \sum_{i=1}^\infty I(b_i) \notag \\
	& = & \sum_{i=1}^\infty (1-h(b_i)) \notag \\
	& = & \sum_{i=1}^\infty (1+b_i \lg{b_i} + (1-b_i)\lg(1-b_i)).
\end{eqnarray}
At any point in time the $b_i$ are only ever specified up to some finite point, after which the $b_i$ are assumed to be completely indeterminate, with $b_i = 0.5$, and hence the information content of all these unspecified digits is nil. The sum (\ref{eqn_info}) is therefore always finite.

This Finite Information Quantity notion of Gisin's is appealing -- as we specify more and more decimal places of any measurement, the digits clearly get more and more uncertain. However, the definition seems like it should bake in the requirement that $|b_i - 0.5| \rightarrow 0$ monotonically, since our certainty about the $(i+1)$st digit would seem to be conditionally dependent on our certainty about the $i$th digit, for all $i$. Moreover, the definition seems to leave unanswered the question of how, and to what accuracy, one is to specify the rational numbers that make up each $b_i$. Finally, while the definition gives an appealing depiction of a limitation on the precision with which a a number can be specified to the right of the decimal point, can we do something similar for very large numbers, in other words regarding the specification of digits to the \textit{left} of the decimal point? Treating digits to the left of the decimal point in the same way in which we have treated digits to the right of the decimal point does not seem to work, but it \textit{does} seem like we can treat the imprecision of a very large integer $N$ simply as the imprecision in $\epsilon$ where $N = \frac{1}{\epsilon}$, with $\epsilon$ as specified in (\ref{eqn_info}).


Returning to Brouwer, and the classical Intuitionistic understanding of mathematics, the notion of time is not just an aspect of the unfolding of digits of real numbers -- the notion extends to all aspects of mathematics.  A theorem starts out with premises and \textit{arrives at} a conclusion -- it is the summary of a mathematical \textit{activity}. It brings two not obviously equal things together, the collection of mathematical objects satisfying the premises and those  satisfying the conclusion, and says they are the same. This activity can only be carried out as long as we have a notion of time -- hence the view of mathematics as activity. Furthermore, from a modern perspective, a theorem is in part information. Information does not exist in the abstract. It must be materialized as bits somewhere -- in the head of a mathematician, on a written page, or stored on a computer, and that materialization as bits requires time and space.

Intuitionism is actually a deep and reasonably well studied subject. There is a conception of natural numbers and, as we have noted, of real numbers, though the conception of each is considerably different from the conventional ones. All infinite collections are considered ``potentially infinite'' and special care has to be taken when reasoning about \textit{all} natural numbers or \textit{all} real numbers. On the other hand, one can reason about \textit{any} natural number or \textit{any} real number. There is a principle of induction for the natural numbers, but a statement about the reals of the form ``every bounded subset of reals has a least upper bound (in $\mathbb{R}$)'' runs into a problem.

Although a further exploration of Intuitionism is beyond the scope of this essay, many questions seem worth delving into: if we rescrutinize the foundations of mathematics and computer science under the assumption that arbitrarily small and arbitrarily large quantities cannot exist, what happens? What parts of the edifice fall apart and what parts can be salvaged? Does Intuitionism form a satisfactory foundation for a least practical computational reasoning? What parts of classical complexity theory does Intuitionism \textit{not} support (if any)? Can we develop complexity theory using Gisin's notion of a Finite Information Quantity, though flipped, as we have indicated, so the notion can apply to large integers?
Finally, are these questions philosophical ones that will never impact the practice of science (and so, perhaps, best be speculated upon at a local pub), or are they, rather, deep scientific questions that could impact our understanding of the universe and even impact the role science has in informing the practice of applied science and engineering? 
I shall delve into these matters with a bit more precision and the specification of three questions for further research at the very end of this essay.

\section{Ultrafinitism}

The notion that there may be a largest natural number has been argued in other ways by various representatives of a school of mathematics known as Ultrafinitism \cite{wiki:ultrafinitism}. Like Intuitionists, Ultrafinitists do not accept the notion of a completed infinite set. Some Ultrafinitists do not accept natural numbers that cannot be expressed in a binary or other decimal expansion that would require more storage space than the assumed storage space of the universe. Others take issue with certain explicit symbolic expressions of natural numbers that they view as potentially uncomputable e.g., taking the floor ($\lfloor x \rfloor$) of a potential upper bound for the first Skewes Number that is dependent on the Riemann Hypothesis: \cite{wiki:skewes} 
\begin{equation*}
	e^{e^{e^{79}}}.
\end{equation*}
Skewes numbers come up in number theory because they are the crossing points of the prime counting number $\pi(x)$, the function giving the number of primes less than a given number $x$, and the logarithmic integral function, $\textrm{li}(x) = \int_0^x \frac{dt}{\ln{x}}$, which, by the Prime Number Theorem, approximates the same number.

Closer to the spirit of this essay is the attempt by Kornai \cite{KORNAI:2003} to characterize the set of numbers expressible via the extraordinarily fast growing Ackermann functions in a theoretical universe capable of storing $2^{512}$ bytes.  $2^{512}$ bytes $\approx10^{155}$ bits, is somewhat more than the   computational capacity of the universe, estimated at $10^{120}$ bits in \cite{LLOYD:2002}. 

A compelling view of Ultrafinitism has been expressed by the Rutgers mathematician Doron Zeilberger \cite{zeilberger:2019}, although he has not written a great deal on the subject. Of note is Zeilberger's observation that integral calculus can be developed perfectly well without limits, and in fact was developed in a much more informal way by its creators Newton and Leibniz in the late 1600s. The epsilon-delta formulation that mathematicians of today are familiar with was not developed until the early 1800s by Balzono and Cauchy \cite{FELSCHER:2000}. In Zeilberger's view a differential equation is just the meaningless ``degenerate case'' of a finite difference equation.

An important contribution to Ultrafinitism is the work on so-called Predicative Arithmetic due to Princeton mathematician Ed Nelson \cite{NELSON:1986}. Nelson found a certain incompatibility in the Peano axioms, with the unconstrained axiom schema of induction and the essentially constructivist axioms defining the properties of the successor function. In his Predicative Arithmetic inductive arguments cannot be ``predicated'' upon numbers that have not yet been shown to exist.
Thus the notion of a natural number being \textit{even} can be defined perfectly well as a natural number $n$ for which there is a previously constructed natural number $m$ and $n=2m$. However, claiming, say, that there is a number divisible by every prime between $1$ and $n$, at stage $n$, is not valid \cite{WILKIE:1990}. It is conceivable that the notion of predicativty has some connection to the P=NP problem (e.g., in distinguishing the languages in P from those in NP).

\section{Concluding Thoughts}

As I end this essay I think it is appropriate to offer 
two quotes about time from Brouwer and Einstein, although the quote by the former will sound dated. I will then wrap up with three of the most pressing questions for future research that arise from our discussion.
As I've noted, time is central to Brouwer's notion of Intuitionism -- a theorem does not make sense outside of time. Brouwer's Intuitionism broke sharply with Hilbert's formalistic school of mathematics \cite{wiki:formalism} (in fact the two men never got along, and had a much publicized falling out in the 1920s over editorial policies at the emminent journal, \textit{Mathematische Annalen}). Hilbert considered mathematics and logic to be statements about the consequences of the manipulation of  sequences of symbols using established manipulation rules. His view was that the symbols could be stripped of any ``interpretation'' and be dealt with purely syntactically.  His aim was to develop a complete and provably consistent axiomatization of all of mathematics. He was the first to articulate how to attempt to do so using a formal language. According to Hilbert, the language must include five components \cite{wiki:formalism}:
\begin{itemize}
\item It must include variables such as x, which can stand for some number. 
\item It must have quantifiers such as the symbol for the existence of an object.
\item It must include equality.
\item It must include connectives such as $\leftrightarrow$ for ``if and only if.''
\item It must include certain undefined terms called parameters. For geometry, these undefined terms might be something like a point or a line, which we still choose symbols for.
\end{itemize}
Of course G\"{o}del's Incompleteness Theorem(s) dealt a death blow to Hilbert's Program. Although the original goal of establishing completeness and consistency has been abandoned, the vast majority of contemporary mathematics still follows Hilbert's paradigm.

\smallskip

Brouwer, looked at mathematics very differently. He did not view axiomatization or even the language of mathematics as paramount. To Brouwer, axiomatization and language, as Hilbert conceived them, were attempts to introduce rigor into an inherently intuitive process. Brouwer declared two founding conceptual ``Acts of Intuitionism'' that gave rise to his various theories. The first of these acts was the following:

\smallskip

\begin{adjustwidth}{2.5em}{2.5em}
``\textit{The First Act of Intuitionism.} Completely separating mathematics from mathematical language and hence from the phenomena of language described by theoretical logic, recognizing that intuitionistic mathematics is an essentially languageless activity of the mind having its origin in the perception of a move of time. This perception of a move of time may be described as the falling apart of a life moment into two distinct things, one of which gives way to the other, but is retained by memory. If the twoity thus born is divested of all quality, it passes into the empty form of the common substratum of all twoities. And it is this common substratum, this empty form, which is the basic intuition of mathematics.'' -- L.E.J. Brouwer \cite[pp.~4--5]{BROUWER:1981}
\end{adjustwidth}

\medskip

Time also plays a central role in the special and general theories of relativity on account of the Lorentzian space-time continuum and the warping of space-time that is the manifestation of the gravitational force. 
The following quote is from a letter Einstein wrote shortly before his death in 1955, upon the passing of his life-long friend and scientific sounding board, Swiss/Italian engineer, Michele Besso. The letter was to Besso's family.

\smallskip

\begin{adjustwidth}{2.5em}{2.5em}
``People like us, who believe in physics, know that the distinction between past, present and future is only a stubbornly persistent illusion.'' -- Albert Einstein \cite{Einstein-Besso}
\end{adjustwidth}

\medskip

 Although, if we take the Intuitionists' word for it and agree that mathematics is inescapably a mental activity intertwined with time, we are brought back, as Einstein was, to the question of `what is time?' In light of the relativity of simultaneity there is no well defined present moment, and to someone traveling at the speed of light relative to us, because of time dilation, our lifetimes go by, to them, in the blink of an eye.
 
 \medskip
 
 Let me now conclude with the three promised questions for future research:
 \begin{enumerate}[(1)]
  \item If we make the bit and the symbols $\bm{0}$ and $\bm{1}$ the foundational construct of all of mathematics, replacing the set and the associated symnbols $\bm{\{}$ and $\bm{\}}$, how should development of the foundations proceed? Should the process be axiomatic or do we take the extreme stance of Brouwer and dispense with an axiomatic approach? Whether or not we side with Brouwer on axiomatics, how does one recover the notion of a set from that of a bit? Should the notion of a set instead just be part of the logical ``plumbing?'' If we try to take Gisin's approach and think of real numbers as finite information quantities it seems that we need a notion of probability before we have a notion of real or other numbers. How can such an idea be made to work `ground-up?' If we consider making the qubit the foundational concept rather than the bit we run into the same problem.
  \item Are Cantor's diagonalization argument, and complexity questions such as P=NP, really meaningful? Note that our ability to separate complexity classes almost completely relies on diagonalization arguments. There is obviously no question that the classical asymptotic analysis of our most common algorithms, say for sorting or finding shortest paths, \textit{is} meaningful and immensely practical. If we are to assert that P=NP is \textit{not} meaningful, we need to do so very precisely and in such a way as to not throw out the baby with the bath water.
 \item In our discussion of inexpressivity, can we be completely precise about the smallest number inexpressible using $n$ bits? After all, there are theorems about such inexpressivity in finite (as well as infinitary) model theory. For example, reachability\footnote{Also known as $s$-$t$ connectivity, in other words that there is a path between two named vertices.} in undirected graphs is known to be inexpressible in the 1st order logic of graphs but expressible in the monadic 2nd order logic of graphs \cite{LIBKIN:2004}.  If we specify a language $L_0$, say the first or second order language of Peano Arithmetic, then the smallest number inexpressible in $n$ bits in $L_0$ will be expressible in $n$ bits in a more expressible language $L_1$ -- a language that can express expressibility in $L_0$ (though not its own internal expressibility/inexpressibility). We thus have $L_1 \supset L_0$. We could seemingly continue in this fashion, obtaining languages $\cdot\cdot\cdot \supset L_2 \supset L_1 \supset L_0$. The fact remains, however, that there is only finitely much we can express with $n$ bits and so at some point we cannot continue the process of taking language extensions. At what point is there a breakdown and how do we state this conundrum formally?
 
 A natural argument to try to make an end run around these inexpressivity limitations is to incorporate time in addition to space. In other words, at time $t=0$ we define a largest (or extremely large) natural number $N_0$ using most, or all, of our storage. Then at time $t=1$ we give a name to $N_0$, and reuse our storage to define a much larger number $N_1$, and so on, obtaining ever larger numbers in time. Since we are assuming finite space, there is then the obvious objection about finite time. However, there is a more serious objection as well. Suppose we have created such an $N_k$ for $k \geq 1$. Since we have destroyed all record of how the $N_0,...,N_{k-1}$ were created our $N_k$ is in fact meaningless.
 \end{enumerate}
 

\bibliographystyle{abbrv}
\bibliography{problem_at_the_foundations_v1.1}

\bigskip

\noindent ``The (Euclidean) point is not the point, but rather, that beyond some point there is no point.'' \\ 
$~~~~~~~~~~~~~~~~~~~~~~~~~~~~~~~~~~~~~~~~~~~~~~~~~~~~~~~~~~~~~~~~~~~~~~~~~~~~~~~~~~~~~~~~~~~~~~~~~~~~~~~~~~~~~~~~~~~~~~~~~~~~$ --JL 

\end{document}